\documentclass[12pt]{article}
\usepackage{amsmath,amssymb,amsfonts,amsthm,graphicx}
\title{\bf A Novel Method for Drawing a Circle Tangent to
Three Circles Lying on a Plane by Straightedge,Compass,and
Inversion Circles}\vspace{2cm}\input{epsf}

\begin{document}
%\vspace{-9cm}
\author{Ahmad Sabihi\thanks{Permanent Address:~First floor,
Bldg.name: Anita,~Golbarg blind alley,~Ershad alley (No.33),~Avecina St.,
Isfahan,~Iran,~Tel:+98(31)34486478,~Cell Phone:+989131115784,~E-mail:~sabihi2000@yahoo.com}\\
\small Professor and researcher at some of Iranian Universities}
\date{}
\maketitle \thispagestyle{empty}
\vspace{-.5cm}
%\textbf{\textit{Dedicated to Professor H.M.Srivastava on the
%%Seventieth Birth Anniversary}}\\\\
{\footnotesize{\bf
Abstract}~~In this paper, we present a novel method to draw
a circle tangent to three given circles lying on a plane.
Using the analytic geometry and inversion (reflection) theorems,
the center and radius of the inversion circle are obtained. Inside
any one of the three given circles, a circle of the similar
radius and concentric with its own corresponding original circle is
drawn.The tangent circle to these three similar circles is
obtained. Then the inverted circles of the three similar circles
and the tangent circle regarding an obtainable point and a
computable power of inversion (reflection) constant are obtained.
These circles (three inverted circles and an inverted
tangent circle)will be tangent together.Just,we obtain another reflection point and power of inversion so that those three reflected circles (inversions of three similar
circles) can be reflections of three original circles, respectively. In such a case,the reflected circle tangent to three reflected circles regarding same new inversion system will be tangent to the three original ones. This circle is our desirable circle. A drawing algorithm is also given for
drawing desirable circle by straightedge and compass. A survey of conformal mapping theory and inversion in higher dimensions is also accomplished. Although, Laguerre transformation might be used for solution of this problem, but we do not make use of this method.  Our novelty is just for drawing a circle tangent to three given circles applying a tangent circle to three identical circles concentric with three given ones and then inverting them as original ones by compass and straightedge not any thing else. \\\\
\textbf{Keywords}: Inversion Circles;Conformal
Mapping;Transformation Geometry;Three-Circle Tangency
Problem;Inversion in Higher Dimension; Laguerre transformation
\section{Introduction}
  ~~~~~The most important and interest problem,which was posed as a challenge
by Apollonius [1] in the third century BC is three-circle tangency
problem.This problem is stated as "given three circles lying on a
certain plane in general,construct all circles tangent to them ".
Apollonius generated the construction for nine cases,briefly,
given three objects, each of them being a circle, a line, or a
point, construct a circle tangent to all three given objects. The
simplest case is the elementary constructs of a circle through
three given points. The most challenging case is that one can draw
a circle tangent to three given circles.That case is of the most
emphasis here. Euclid and Archmedes had presented the construction
for three points and three lines,but Apollonius gave the
construction methods in his own books I and II on three-circle
problem as well. Altshiller-Court [2]in 1925 gave a construction of
same problem by reducing it to the construction of PCC (Point
Circle-Circle). Fran\c{c}ois Vi\`{e}te, Isaac Newton, and
Joseph-Diaz Gergonne [3] individually solved the problem by
different approaches.Vi\`{e}te made use of the geometric
fundamentals, Newton the conic approach, and Gergonne developed a
new understanding of inversion geometry properties.Vi\`{e}te's
work is fine for any three given circles,but not in an environment,
where the given conditions can change after completion of the
construction. Newton's solution is also a bit problematic when it
comes to dynamic geometry software. A change in the size or
location of one of the given circles is likely to capsize the
construction. His construction method is not unique so that
decision must be made along the construction way. The method is
based on the intersection of hyperbolas. An intersection of two
hyperbolas is the case not three. Two hyperbolas can intersect at
four points, and only two of which solve problem [3]. Gergonne in
1816 [3] published an inversion-based solution to the tangency
problem. The Gergonne's solution can be applied to many of the
nine other cases.His method is the most exhaustive solution,which
is specially flexible regarding positions of the given
circles.\\\indent Inversion was invented by J.Steiner around 1830
and transformation by inversion in circle was invented by
L.I.Magnus in 1831 [5]. For further information, refer to [6-9].
The author,surveys on inversive geometry,transformation theory,
and conformal mapping in the Section 2. A novel geometrical method
for drawing a circle tangent to the three given circles is given
in the Section 3. A new algorithm for drawing same by Compass,and
Straightedge without using any other axillary tools comes in the
Section 4.\\
\setcounter{equation}{0}
\section{Inversive Geometry,Transformation Theory,and Conformal Mapping}
 ~~~~~Let $C_{1}$ be a circle passing through the point $A$
 perpendicular to $C_{2}$. Let $\acute{A}$ be a second point of
 intersection of $OA$ with $C_{1}$,then $OC$ is tangent to
 $C_{1}$. This means that $OA.O\acute{A}=OC^{2}=R^{2}$ where $R$ is
 the radius of $C_{2}$. $C_{2}$ and $O$ are called the reflection
 circle and center respectively. Also,the location of point
 $\acute{A}$ does not depend on circle $C_{1}$,because $\acute{A}$
 lies on the line $OA$ at distance equal to $\frac{R^{2}}{OA}$. The
 point $\acute{A}$ is known as the
 \textit{inverse},\textit{inversion},or \textit{reflection} of $A$
 with respect to $C_{2}$. Therefore,$\acute{A}$ is the reverse or
 reflect of $A$ and vice versa. Except the center of $C_{2}$,which
 is called as the center of inversion,all points as $A$ in the
 plane have inverse images as $\acute{A}$. The term inversion also
 applies to define transformation of the plane.This transformation
 maps the points inside $C_{2}$ to the points in its exterior and
 vice versa. The center of inversion is often [4-5] left over as a
 point with no inverse image,but sometimes is said to be mapped to
 the point at infinity. The origin or zero in
the circle inversion mapping requires a special definition to
adjoin a point at infinity as $\infty$. $R$ is called the radius
of inversion,$R^{2}$ as the power of inversion. As stated in the
introduction,inversive geometry refers to a study of
transformation geometry by Euclidean transformations together with
inversion in an n-sphere so:
\begin{equation}
x_{i}\mapsto \frac{r^{2}x_{i}}{\sum_{j=1}^{n}x_{j}^{2}}\
\end{equation}
where $r$ denotes the radius of the inversion. In two
dimensions,with $r=1$,this is indeed a circle inversion with
respect to the unit circle so that $x\mapsto
\frac{r^{2}x}{|x|^{2}}$. If two inversions in concentric circles
are combined,then a similarity,homothetic transformation will be
resulted by the ratio of the circle radii:
\begin{equation}
x\mapsto \frac{r^{2}x}{|x|^{2}}=y\mapsto
\frac{k^{2}y}{|y|^{2}}=(\frac{k}{r})^{2}x\
\end{equation}
The algebraic form of the inversion in a unit circle is
$w=\frac{1}{\bar{z}}$ because
\begin{equation}
\frac{\bar{z}}{|z|}.\frac{z}{|z|}=R^{2}=1\
\end{equation}
and the reciprocal of $Z$ is $\frac{\bar{z}}{z}$. In the complex
plane,reciprocation leads to the complex projective line and is
often called the Riemann sphere. The subspaces and sub groups of
this space and group of mappings,were applied to produce early
models of hyperbolic geometry given by Arthur Cayley,Felix Klein,
Henri Poincar\'{e}. Therefore,inverse geometry includes the ideas
originated by Lobachevsky and Bolyai in their plane geometry. More
generally, a map $f:U\longrightarrow V$ is called conformal or
angle-preserving at $v_{0}$,if it preserves oriented angles
between curves through $v_{0}$ with respect to their
orientations. Conformal maps preserve both angles and shapes of
infinitesimally small figures,but not necessarily their size. Some
examples of conformal maps could be made from complex analysis. Let
$A$ be an open subset of the complex plane,$\Bbb C$,then a
function $g:A\longrightarrow \Bbb C$ is defined as a conformal map
if and only if it is holomorphic and its derivation is everywhere
non-zero on $A$.But if $g$ is anti-holomorphic (conjugate to a
holomorphic function),it still preserves angles except it reverses
their orientation. Conformal mapping application in electromagnetic
potential,heat conduction,and flow of fluids are taken into
considerations. The ideas are that the problems at hand with a
certain geometry could be mapped into a problem with simpler
geometry or a geometry which have already been solved. A such
mappings are found by transformation [10]. For further study,refer
to [11-12]. Laguerre transformation [13] might be used for solving our problem using a plane spanned by three points in a three dimensional space and transforming it into a horizontal plane, but we do not make use of such transformations. Also, methods of Cyclography might be applied to solve it.  Our novelty is just for drawing a circle tangent to three given circles applying a tangent circle to three identical circles concentric with three given ones and then inverting them as original ones by compass and straightedge not any thing else. 
\setcounter{equation}{0}
\section{The Novel Method}
\textbf{3.1 General Descriptions}\\\\
\indent Let $C_{1},C_{2}$,and $C_{3}$ be three given circles of
optional radii $R_{1}<R_{2}<R_{3}$. Let's draw three similar circles of radii $R$
concentric with each one of circles so that inequality $R<R_{1}<R_{2}<R_{3}$ holds. Adjoin the
centers of circles together by straight-lines. Then draw
perpendicular bisector of the created line-segments. The intersection point of
perpendicular bisectors is indeed the center of tangent circle to
the three similar circles of the radii $R$ since it is
equidistance from the circles' centers. Also, having equal
radii,the distances of the intersection point from three circles
are identical. This means that we able to draw a tangent
circle to the three similar circles from the intersection point. Just, let the problem be
solved. Using the reflection theorems and methods,the reflected
objects (or reflected circles) of three similar circles of radii
$R$ with respect to a fixed point (where will be obtained later by
drawing) are obtained. In this paper,we carry out following steps 3.1.1 through 3.1.7  based on  Fig.1\\\\
\textbf{3.1.1}Consider the three given circles $C_{1},C_{2}$,and $C_{3}$ of radii $R_{1},R_{2}$ and ,$R_{3}$, respectively.\\\\
\textbf{3.1.2}Draw the concentric circles of radii $R$ denoted by
$\underline{C_{1}},\underline{C_{2}}$,and $\underline{C_{3}}$
(these symbols are not shown in Fig.1) inside the circles
$C_{1},C_{2}$,and $C_{3}$, respectively.\\\\
\textbf{3.1.3}Draw the tangent circle $C_{4}$ to the three similar circles of the radii $R$, if their centers are not on a line\\\\
\textbf{3.1.4}Draw the reflected circles of the three similar circles of radii $R$ denoted by $\acute{C}_{1},\acute{C}_{2}$,and $\acute{C}_{3}$ with respect to the reflection center $O$\\\\
\textbf{3.1.5}Draw the reflected circle of the tangent circle $C_{4}$ i.e.$\acute{C}_{4}$. This reflected circle is also tangent to the reflected circles of the three similar circles of the radii $R$.\\\\
\textbf{3.1.6}Define a new center and a power of inversion so that the three reflected circles $\acute{C}_{1},\acute{C}_{2}$,and $\acute{C}_{3}$ (created by three similar circles of radii $R$) can be as the new reflected circles of $C_{1},C_{2}$,and $C_{3}$. The new reflection center is identical to the old one.\\\\
\textbf{3.1.7}Draw the reflection of the tangent circle $\acute{C}_{4}$ with respect to the new given center (or same old one) and new reflection constant (power of inversion value). This circle denotes ${C''}_{4}$ and is a desirable circle in this paper, where is tangent to the three given circles.\\\\
\textbf{3.1.8}If the centers of the three given circles (and three
similar circles)$C_{1},C_{2}$,and $C_{3}$ are on a straight line,then we
should draw a tangent line to the three similar circles instead of
a tangent circle and obtain the reflected circle of this line and
continue steps 3.1.5 to 3.1.7.\\\\
\textbf{3.2 Determination of the Center and Radius of Reflection Circle}\\\\
\begin{figure}
  \centering
  \includegraphics[width=16cm]{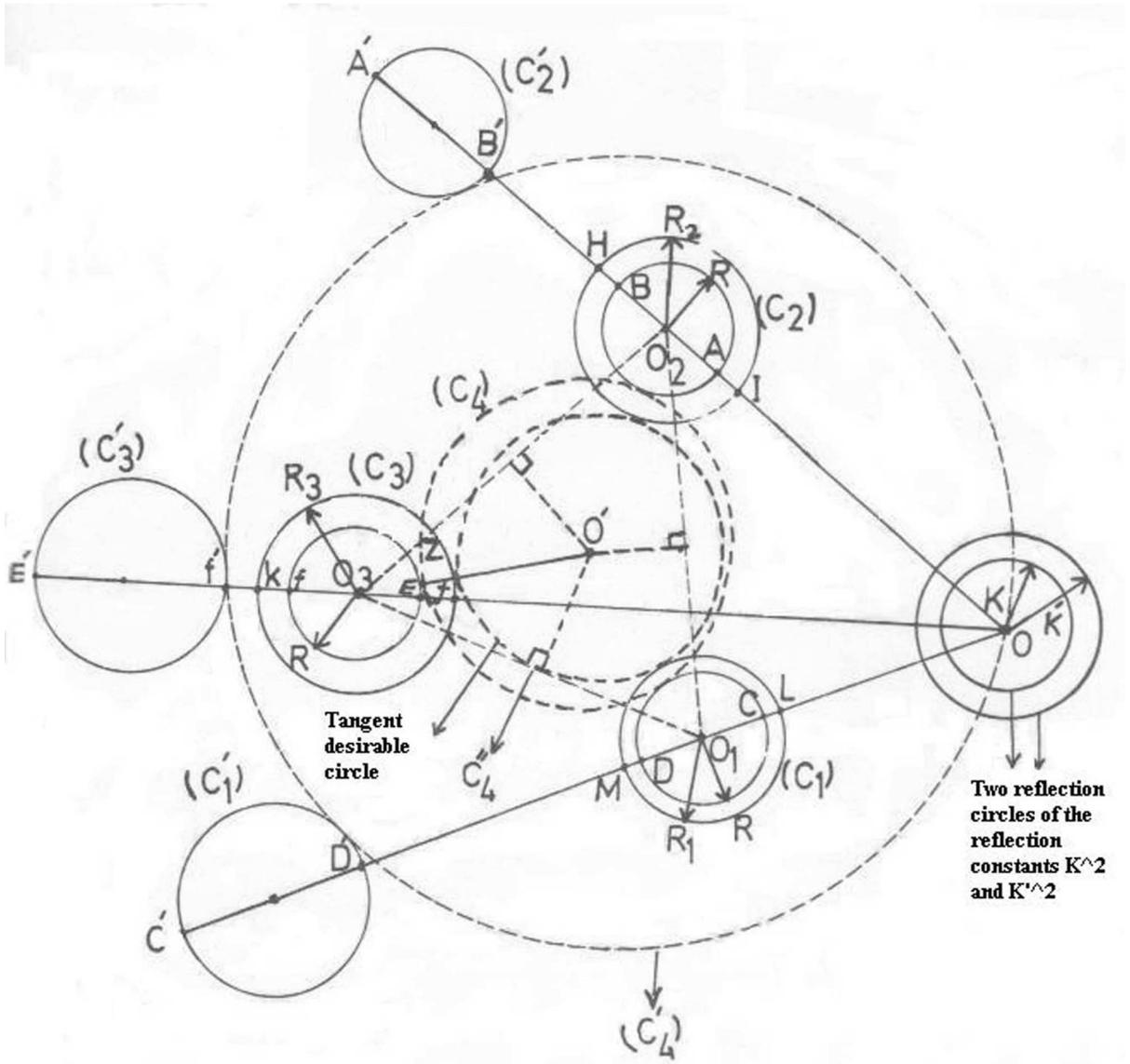}\\
  \caption{Original, identical
circles, their inversions, and two reflection circles of the reflection constants (or power of inversion) $K^{2}$ and
$\acute{K}^{2}$}\label{Fig:c1}
\end{figure}
\textbf{Theorem1,}
Let we have the first reflection center $O$ referring Fig.1. We know
the following relations hold:\\
\begin{equation}
OA.O\acute{A}=OB.O\acute{B}=OC.O\acute{C}=OD.O\acute{D}=OE.O\acute{E}=Of.O\acute{f}=k^{2}
\end{equation}\ where $k^{2}$ denotes the power of inversion of the three similar
circles. The points
$\acute{A},\acute{B},\acute{C},\acute{D},\acute{E}$ and
$\acute{f}$ denote reflected points for points $A$,$B$,$C$,$D$,$E$ and $f$ on the three reflection circles $\acute{C_{1}},\acute{C_{2}}$,and $\acute{C_{3}}$, respectively. Then, to make
the reflected circles $\acute{C_{1}},\acute{C_{2}}$,and
$\acute{C_{3}}$ as the reflections of the circles $C_{1},C_{2}$,
and $C_{3}$, respectively, a new center of reflection can be defined so that
lies on the previous one. Hence,the center $O$ should be the
reflection center of both reflection systems so that:\\\\
\begin{equation}
OI.O\acute{A}=OH.O\acute{B}=OL.O\acute{C}=OM.O\acute{D}=OJ.O\acute{E}=OK.O\acute{f}=\acute{k}^{2}
\end{equation}\
\textbf{\textit{Proof}}\\\\
To have $O$ as the common center of both reflection systems:
\begin{equation}
OI.O\acute{A}.OH.O\acute{B}=\acute{k}^{4}\
\end{equation}\
\begin{equation}
OI.OH.\frac{k^{2}}{OA}.\frac{k^{2}}{OB}=\acute{k}^{4}\ \Rightarrow
\frac{OH.OI}{OA.OB}=\frac{OL.OM}{OC.OD}=\frac{OJ.OK}{OE.Of}=\frac{\acute{k}^{4}}{{k}^{4}}=M^{2}\
\end{equation}\
where $M^{2}$ denotes the ratio of the second power of the
inversions of both reflection systems. This is a constant value.
It also denotes the ratio of the powers of the point $O$ with
respect to concentric circles of radii $(R_{1},R),(R_{2},R)$,and
$(R_{3},R)$. This means that ratio of the tangent line segments
drawn from the point $O$ over both concentric circles $(R_{1},R),(R_{2},R)$,and
$(R_{3},R)$ should be identical. The problem is that we wish to
find geometrical location of all points whose powers are
identical with respect to two given concentric circles. Referring
Fig. 1 and Fig.2,we have the relations (5-3).\\\\\\
\begin{figure}
  \centering
  \includegraphics[width=8cm]{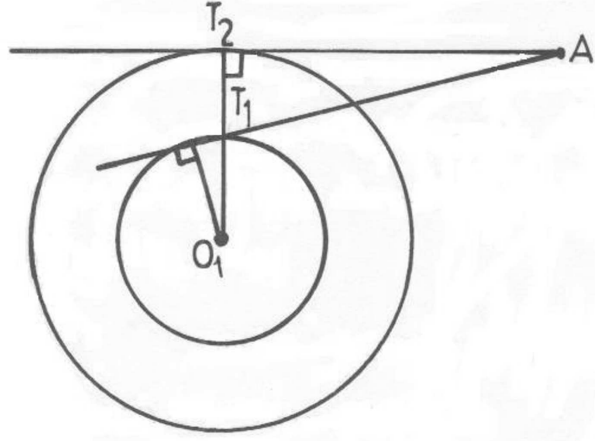}\\
  \caption{Geometrical location of all points whose point's powers are identical with respect to
two given concentric circles} \label{Fig:c2}
\end{figure}

\begin{equation}
\frac{OH.OI}{OA.OB}=(\frac{AT_{2_2}}{AT_{1_2}})^{2},~~\frac{OL.OM}{OC.OD}=(\frac{AT_{2_1}}{AT_{1_1}})^{2}~~\frac{OJ.OK}{OE.Of}=(\frac{AT_{2_3}}{AT_{1_3}})^{2}\
\end{equation}\
Therefore,using (4-3) and (5-3)
\begin{equation}
(\frac{AT_{2_1}}{AT_{1_1}})^{2}=(\frac{AT_{2_2}}{AT_{1_2}})^{2}=(\frac{AT_{2_3}}{AT_{1_3}})^{2}=M^{2}\
\end{equation}\
where $AT_{1_1}$ to $AT_{1_3}$ denote tangent lines drawn from the
point $A$ or same $O$ (a typical reflection center) to the three
similar circles of the radii $R$,respectively. Also, $AT_{2_1}$ to $AT_{2_3}$
denote tangent lines from same point to the circles
$C_{1},C_{2}$,and $C_{3}$, respectively. For example, for circles
$C_{1}$ and $\underline{C_{1}}$ (two concentric circles shown in Fig.2) of the
center $O_{1}$ we have:
\begin{equation}
\frac{OL.OM}{OC.OD}=(\frac{AT_{2_1}}{AT_{1_1}})^{2}=M^{2}=\frac{AO_{1}^{2}-R_{1}^{2}}{AO_{1}^{2}-R^{2}}\Rightarrow
AO_{1}^{2}-R_{1}^{2}=M^{2}.AO_{1}^{2}-M^{2}.R^{2}\Rightarrow
AO_{1}=\sqrt{\frac{R_{1}^{2}-M^{2}.R^{2}}{1-M^{2}}}\
\end{equation}\
The geometrical location of the point $A$ is a circle of the
radius $AO_{1}$ and the center $O_{1}$. For the two other
concentric circles $(C_{2},\underline{C_{2}})$,and
$(C_{3},\underline{C_{3}})$, the method is similar. Due to
$0<M^{2}<1$ we can choose an optional value e.g.
$M^{2}=\frac{1}{2}$. We know $O_{2}$ and $O_{3}$ are the centers and
$R_{2}$ and $R_{3}$, the radii of two other concentric circles,
respectively, therefore regarding (7.3), we have 
\begin{equation}
AO_{1}^{2}-2R_{1}^{2}=-R^{2},~~AO_{2}^{2}-2R_{2}^{2}=-R^{2},~~AO_{3}^{2}-2R_{3}^{2}=-R^{2}\
\end{equation}\
Let the representing parametric equations of $C_{1},C_{2}$,and
$C_{3}$ as follows:
\begin{equation}
C_{1}(x,y):~~x^{2}+y^{2}+a_{1}x+b_{1}y+c_{1}=0,~~C_{2}(x,y):~~x^{2}+y^{2}+a_{2}x+b_{2}y+c_{2}=0,~~C_{3}(x,y):~~x^{2}+y^{2}+a_{3}x+b_{3}y+c_{3}=0\
\end{equation}\
Let $A(x_{0},y_{0})$ be the mentioned typical reflection center. Using
the equations (9-3),we obtain the circle centers 
$O_{1}(\frac{-a_{1}}{2},\frac{-b_{1}}{2})$,$O_{2}(\frac{-a_{2}}{2},\frac{-b_{2}}{2})$,
and $O_{3}(\frac{-a_{3}}{2},\frac{-b_{3}}{2})$ then
\begin{equation}
AO_{1}^{2}=(x_{0}+\frac{a_{1}}{2})^{2}+(y_{0}+\frac{b_{1}}{2})^{2},~~R_{1}^{2}=\frac{1}{4}(a_{1}^{2}+b_{1}^{2})-c_{1}\Rightarrow
AO_{1}^{2}-2R_{1}^{2}=x_{0}^{2}+y_{0}^{2}+a_{1}x_{0}+b_{1}y_{0}+c_{1}+(c_{1}-\frac{1}{4}(a_{1}^{2}+b_{1}^{2}))\
\end{equation}\
By the same method,we obtain
\begin{equation}
AO_{2}^{2}-2R_{2}^{2}=x_{0}^{2}+y_{0}^{2}+a_{2}x_{0}+b_{2}y_{0}+c_{2}+(c_{2}-\frac{1}{4}(a_{2}^{2}+b_{2}^{2}))\
\end{equation}\
\begin{equation}
AO_{3}^{2}-2R_{3}^{2}=x_{0}^{2}+y_{0}^{2}+a_{3}x_{0}+b_{3}y_{0}+c_{3}+(c_{3}-\frac{1}{4}(a_{3}^{2}+b_{3}^{2}))\
\end{equation}\\\\
\textbf{Theorem2,}
let $K_{1}(x,y)$, $K_{2}(x,y)$ and $K_{3}(x,y)$ denote the new
circle's equations eliminating zero indices as follows, then the center of inversion for three given circles is obtained by intersection of normal lines created by these circles. 
\begin{equation}
K_{1}(x,y)=AO_{1}^{2}-2R_{1}^{2}+R^{2}=C_{1}(x,y)+(c_{1}-\frac{1}{4}(a_{1}^{2}+b_{1}^{2}))+R^{2}=0\
\end{equation}\
\begin{equation}
K_{2}(x,y)=AO_{2}^{2}-2R_{2}^{2}+R^{2}=C_{2}(x,y)+(c_{2}-\frac{1}{4}(a_{2}^{2}+b_{2}^{2}))+R^{2}=0\
\end{equation}\
\begin{equation}
K_{3}(x,y)=AO_{3}^{2}-2R_{3}^{2}+R^{2}=C_{3}(x,y)+(c_{3}-\frac{1}{4}(a_{3}^{2}+b_{3}^{2}))+R^{2}=0\
\end{equation}
\textbf{\textit{Proof}}\\\\
Consider\\
\begin{equation}
K_{3}-K_{2}:~~C_{3}(x,y)-C_{2}(x,y)+c_{3}-c_{2}+\frac{1}{4}(a_{2}^{2}-a_{3}^{2}+b_{2}^{2}-b_{3}^{2})=0\
\end{equation}
\begin{equation}
K_{2}-K_{1}:~~C_{2}(x,y)-C_{1}(x,y)+c_{2}-c_{1}+\frac{1}{4}(a_{1}^{2}-a_{2}^{2}+b_{1}^{2}-b_{2}^{2})=0\
\end{equation}
\begin{equation}
K_{3}-K_{1}:~~C_{3}(x,y)-C_{1}(x,y)+c_{3}-c_{1}+\frac{1}{4}(a_{1}^{2}-a_{3}^{2}+b_{1}^{2}-b_{3}^{2})=0\
\end{equation}
The equations (16-3) to (18-3) are indeed the equations of the
normal lines to the lines which adjoin the centers of both circles
$(C_{3},C_{2})$, $(C_{2},C_{1})$, and $(C_{3},C_{1})$, respectively.
Take the points $M_{1},M_{2}$ and $M_{3}$ placing at the bisector
of the line segments constructed, then referring to
Fig.3, we could make the relations (19-3).\\\\

\begin{figure}
  \centering
  \includegraphics[width=16cm]{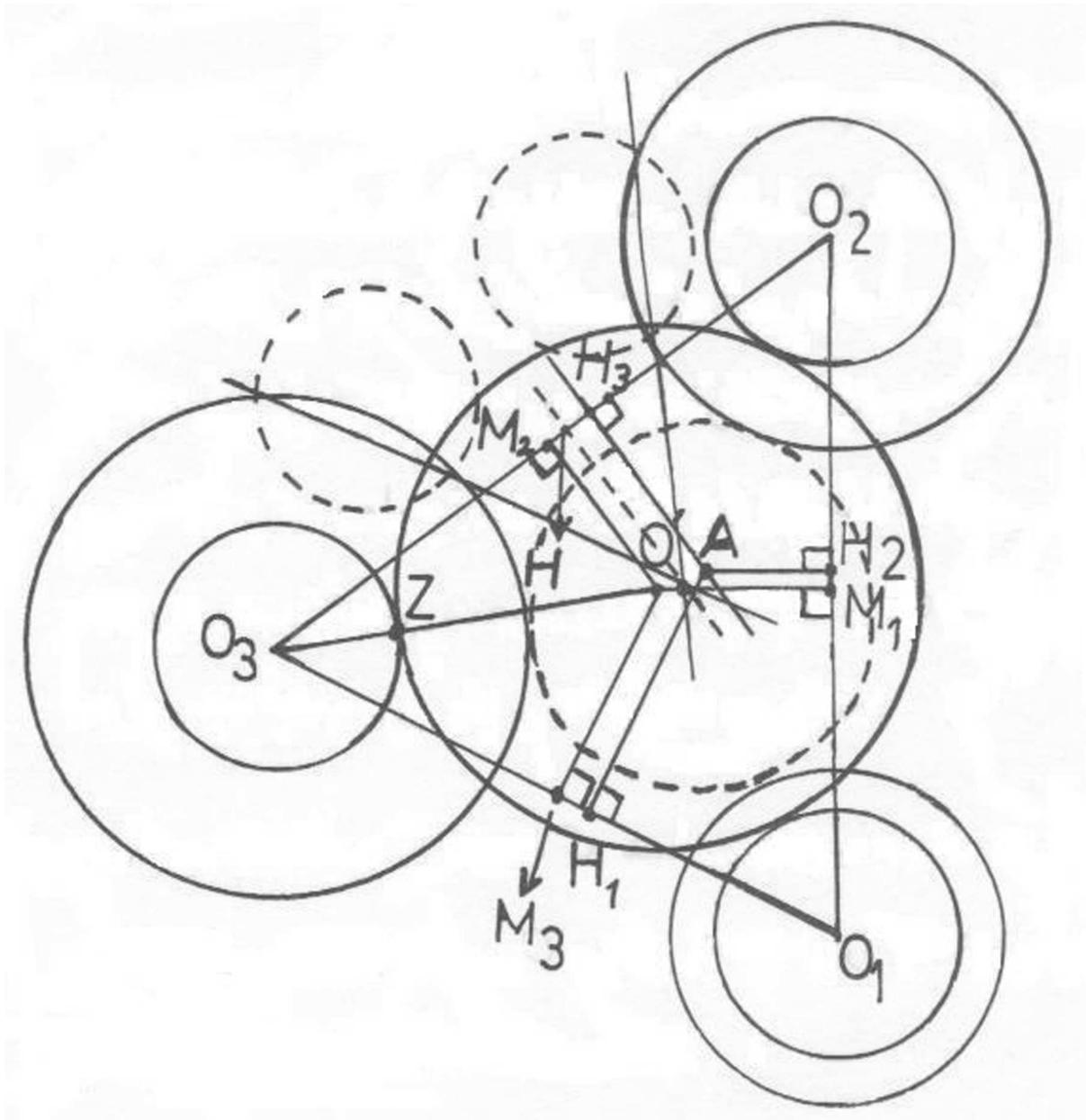}\\
  \caption{The method of
determination of the reflection center $A$} \label{Fig:c3}
\end{figure}

\begin{equation}
AO_{2}^{2}-AO_{1}^{2}=2\overline{O_{2}O_{1}}.\overline{M_{1}H_{2}},~~AO_{3}^{2}-AO_{2}^{2}=2\overline{O_{3}O_{2}}.\overline{M_{2}H_{3}},~~AO_{3}^{2}-AO_{1}^{2}=2\overline{O_{3}O_{1}}.\overline{M_{3}H_{1}}\
\end{equation}\
Therefore, the lines passing through the points $H_{1},H_{2}$,and
$H_{3}$ denote our desirable geometrical locations of the points. These locations point out to the set of the points, whose difference of their power of distances from the centers of each
two circles are constant values. These are the lines perpendicular to the
line segments adjoining their two centers. The distance between
the perpendicular line and the center of the line segment for each two
circles is:
\begin{equation}
\overline{M_{1}H_{2}}=\frac{R_{2}^{2}-R_{1}^{2}}{\overline{O_{2}O_{1}}},~~\overline{M_{2}H_{3}}=\frac{R_{3}^{2}-R_{2}^{2}}{\overline{O_{3}O_{2}}},~~\overline{M_{3}H_{1}}=\frac{R_{3}^{2}-R_{1}^{2}}{\overline{O_{3}O_{1}}}\
\end{equation}\
The given lines presented by the equations (16-3) to (18-3) should
intersect each other at a unique point because
\begin{equation}
(K_{3}-K_{2})+(K_{2}-K_{1})=(K_{3}-K_{1})\
\end{equation}\
Therefore, from the intersection of the two lines $(K_{3}-K_{2})$ and
$(K_{2}-K_{1})$,the line $(K_{3}-K_{1})$ is obtained.Thus,the
point $(x,y)$ satisfies three line equations (16-3)
through (18-3). This means that the center of inversion is found
and the problem is solved. \setcounter{equation}{0}
\section{The Novel Drawing Method by the Straightedge and Compass in form of a New Algorithm}
The algorithm is given regarding Figs.1,3\\\\
\textbf{4.1}Draw three similar concentric circles
$\underline{C_{1}},\underline{C_{2}}$,and $\underline{C_{3}}$ of
radii $R$ inside the three given circles $C_{1},C_{2}$,and $C_{3}$,
respectively so that $R<R_{1}<R_{2}<R_{3}$.\\\\
\textbf{4.2}Adjoin the centers of the circles $C_{1},C_{2}$,and $C_{3}$ to each other.\\\\
\textbf{4.3}Draw the perpendicular bisectors of the line segments
adjoining each of two centers to meet at a point as ${O''}$(not shown in Fig.3).\\\\
\textbf{4.4}Adjoin the point ${O''}$ and center of one
of the circles,say,$\underline{C_{3}}$ (Fig.3) to meet at the point $Z$.\\\\
\textbf{4.5} Draw a circle of the center ${O''}$ and radius ${O''}Z$ to be
tangent to the three similar circles $\underline{C_{1}},\underline{C_{2}}$,and $\underline{C_{3}}$.\\\\
\textbf{4.6} Regarding Fig.3,the reflection center point $A$
could be found as follows:\\\\
Draw the radical axis of both circles
$(C_{1},C_{2})$,$(C_{1},C_{3})$, and $(C_{2},C_{3})$.This can be
done by drawing two optional circles to be intersected with
circles,say,$C_{2}$,and $C_{3}$. These two radical axes 
intersect with one another at the point $\acute{O}$. A perpendicular 
line should be drawn from this point to the line segment
$O_{2}O_{3}$ to meet it at a point. This point is called
$\acute{H}$. Lie the compass's needle on the $\acute{H}$ and open
its other arm by $\acute{H}M_{2}$ based on what is shown in the
Fig.3. Then, draw an arc toward smaller circle side,say,$C_{2}$,
to be intersected with $O_{2}O_{3}$ at the point $H_{3}$. Draw a
perpendicular line to $O_{2}O_{3}$ from the point $H_{3}$. This line is indeed
representing the line $(K_{3}-K_{2})$. The method should be
repeated for drawing the lines $(K_{2}-K_{1})$, and
$(K_{3}-K_{1})$ by same method. These lines 
will intersect each other at the point $A$,which denotes the reflection center.\\\\
\textbf{4.7} Obtain the reflections of three similar circles of
equal radii $R$ (i.e. $\underline{C_{1}},\underline{C_{2}},
\underline{C_{3}}$)and $C_{4}$ with respect to the reflection
center $A$ and power of inversion $K^{2}$. These reflections
denote $\acute{C_{1}},\acute{C_{2}}, \acute{C_{3}}$,and
$\acute{C_{4}}$ according to Fig.1. The circle $\acute{C_{4}}$ is
tangent to the circles $\acute{C_{1}},\acute{C_{2}}$,and
$\acute{C_{3}}$.\\\\
\textbf{4.8}According to given proof, three circles
$\acute{C_{1}},\acute{C_{2}}$,and $\acute{C_{3}}$ can also be as
reflections of the three circles $C_{1},C_{2}$,and $C_{3}$ with regard to same reflection center and by different power of inversion $\acute{K}^{2}$.\\\\
\textbf{4.9}Based on the explanations given in 4.8 and referring
to the relation (4-3), we find out that
$\frac{\acute{k}^{4}}{{k}^{4}}=M^{2}$. Therefore, the new power of
inversion is obtained by the relation $\acute{k}^{2}=M.{k}^{2}$.
Both ${k}^{2}$ and $M$ are known,thus $\acute{k}^{2}$ is also
known, because, we can obtain size $M$ by the relation
$M^{2}=\frac{1}{2}$ and by means of straightedge and compass
tools. This means that the reflections of the three given circles
$\acute{C_{1}},\acute{C_{2}}$,and $\acute{C_{3}}$ of the center
$A$ and constant $\acute{k}^{2}$ can
be the three main circles $C_{1},C_{2}$,and $C_{3}$, respectively.\\\\
\textbf{4.10}Referring to Fig.1,the reflection of the circle
$\acute{C_{4}}$ is ${C''_{4}}$,which is a desirable solution.
${C''_{4}}$ is tangent to the three circles $C_{1},C_{2}$,and
$C_{3}$ due to ${C''_{4}}$ is also reflection of
$\acute{C_{4}}$ with regard to the point
$A$ and the reflection constant $\acute{k}^{2}$.Therefore,the problem of three-circle tangency is completely drawn by compass and straightedge.\\\\
Only tools applying in this method to draw desirable circle,are both Compass and straightedge
without any additional ones. Our method may be extended
to three-dimension space by turning the plane of each one of
circles about their related diameter. In such a case,we can make a
sphere tangent to three given spheres. Also,this method would also
be interested if it could be extended to higher dimension spheres.
If three original circles are of concentric circles or only one of
three circles is placed inside another one and the other one is
outside of both,then we will not be able to draw a tangent circle
to them. Further investigations about drawing limitations are left
for readers.

\end{document}